\documentclass[12pt]{article}
\usepackage{amsmath}
\usepackage{amssymb}
\renewcommand\smallskip{\vskip\smallskipamount}
\renewcommand\medskip{\vskip\medskipamount}
\renewcommand\bigskip{\vskip\bigskipamount}

\begin{document}

\footnotetext{The author is partially supported by NSF Grant
DMS-0707086 and a Sloan Research Fellowship.}

\begin{center}
\begin{large}
\textbf{A Penrose-Like Inequality for General Initial Data Sets}
\end{large}

\bigskip\bigskip
MARCUS A. KHURI
\bigskip\bigskip
\end{center}

\begin{abstract}
We establish a Penrose-Like Inequality for general (not
necessarily time symmetric) initial data sets of the Einstein
equations which satisfy the dominant energy condition.  More
precisely, it is shown that the ADM energy is bounded below by an
expression which is proportional to the square root of the area of
the outermost future (or past) apparent horizon.
\end{abstract}

\bigskip

\begin{center}
\textbf{1.  Introduction}
\end{center}\setcounter{equation}{0}
\setcounter{section}{1}

  In an attempt to find a counterexample for his Cosmic Censorship
Conjecture, R. Penrose [12] proposed a necessary condition for its
validity, in the form of an inequality relating the ADM mass and
area of any event horizon in an asymptotically flat spacetime:
\begin{equation}
\mathrm{Mass}\geq\sqrt{\frac{\mathrm{Area}}{16\pi}}.
\end{equation}
Unfortunately this Penrose Inequality can only be proven with
knowledge of the full spacetime development, as otherwise it would
not be possible to locate the event horizon in a given spacelike
slice.  Thus it is customary to reformulate (1.1) so that the
quantities involved may be calculated solely from local
information, namely initial data sets for the Einstein equations.
By an initial data set we are referring to a triple $(M,g,k)$,
consisting of a Riemannian 3-manifold $M$ with metric $g$ and a
symmetric 2-tensor $k$ representing the extrinsic curvature of a
spacelike slice. These data are required to satisfy the constraint
equations
\begin{eqnarray*}
16\pi\mu&=&R+(\mathrm{Tr}_{g}k)^{2}-|k|^{2},\\
8\pi J_{i}&=&\nabla^{j}(k_{ij}-(\mathrm{Tr}_{g}k)g_{ij}),
\end{eqnarray*}
where $R$ is scalar curvature and $\mu$, $J$ are respectively the
energy and momentum densities for the matter fields.  If all
measured energy densities are nonnegative then $\mu\geq|J|$, which
will be referred to as the dominant energy condition.  Moreover
the initial data set will be taken to be asymptotically flat (with
one end), so that at spatial infinity the metric and extrinsic
curvature satisfy the following fall-off conditions
\begin{equation*}
|\partial^{l}(g_{ij}-\delta_{ij})|=O(r^{-l-1}),\text{ }\text{
}\text{ }\text{ }|\partial^{l}k_{ij}|=O(r^{-l-2}),\text{ }\text{
}\text{ }l=0,1,2,\text{ }\text{ }\text{ as }\text{ }\text{
}r\rightarrow\infty.
\end{equation*}
The ADM energy and momentum are then well defined by
\begin{equation*}
E=\lim_{r\rightarrow\infty}\frac{1}{16\pi}\int_{S_{r}}(\partial_{i}g_{ij}-\partial_{j}g_{ii})\nu^{j},\text{
}\text{ }\text{ }\text{
}P_{i}=\lim_{r\rightarrow\infty}\frac{1}{8\pi}\int_{S_{r}}(k_{ij}-(\mathrm{Tr}_{g}k)g_{ij})\nu^{j},
\end{equation*}
where $S_{r}$ are coordinate spheres in the asymptotic end with
unit outward normal $\nu$.\par
  The strength of the gravitational field in the vicinity of a
2-surface $\Sigma\subset M$ may be measured by the null expansions
\begin{equation*}
\theta_{\pm}:=H_{\Sigma}\pm\mathrm{Tr}_{\Sigma}k,
\end{equation*}
where $H_{\Sigma}$ is the mean curvature with respect to the unit
outward normal (pointing towards spatial infinity). The null
expansions measure the rate of change of area for a shell of light
emitted by the surface in the outward future direction
($\theta_{+}$), and outward past direction ($\theta_{-}$).  Thus
the gravitational field is interpreted as being strong near
$\Sigma$ if $\theta_{+}< 0$ or $\theta_{-}< 0$, in which case
$\Sigma$ is referred to as a future (past) trapped surface.
Future (past) apparent horizons arise as boundaries of future
(past) trapped regions and satisfy the equation $\theta_{+}=0$
($\theta_{-}=0$). In the setting of the initial data set
formulation of the Penrose Inequality, apparent horizons take the
place of event horizons, in that the area of the event horizon is
replaced by the area of the outermost apparent horizon (or in some
formulations by the least area required to enclose an apparent
horizon).\par
  The Penrose Inequality has been established by Huisken \& Ilmanen
[8] and by Bray [2] in the time symmetric case, that is when
$k=0$.  At the present time the conjecture for arbitrary initial
data sets remains open, however recently Bray and the author [3]
have succeeded in reducing this problem to the question of
existence for a canonical system of partial differential
equations.  The purpose here is to establish a Penrose-Like
Inequality for arbitrary initial data satisfying the dominant
energy condition. This new inequality will generalize the
following one obtained by Herzlich [7] in the time symmetric
case.\medskip

\textbf{Theorem 1.1.}  \textit{Let $(M,g)$ be a 3-dimensional
asymptotically flat Riemannian manifold with nonnegative scalar
curvature, and boundary consisting of a minimal 2-sphere with area
$|\partial M|$.  Then the ADM energy satisfies}
\begin{equation*}
\mathrm{E}(g)\geq\frac{\sigma}{2(1+\sigma)}\sqrt{\frac{|\partial
M|}{\pi}}
\end{equation*}
\textit{where}
\begin{equation*}
\sigma=\sqrt{\frac{|\partial M|}{\pi}}\inf_{v\in
C_{c}^{\infty}\atop v \neq 0}\frac{\parallel\nabla
v\parallel^{2}_{L^{2}(M)}}{\parallel
v\parallel^{2}_{L^{2}(\partial M)}}.
\end{equation*}
\textit{Furthermore equality holds if and only if $(M,g)$ is a
portion of the $t=0$ slice of the Schwarzchild spacetime with mass
$\sqrt{|\partial M|/16\pi}$.}\medskip

%\textbf{Remark 1.2.}  \textit{For simplicity, all asymptotically
%flat manifolds in these notes are assumed to have just one
%end.}\medskip

  A useful device for extending results in the time symmetric case
to the general case is Jang's deformation [9] of the initial data,
which was successfully employed by Schoen and Yau [13] in their
proof of the Positive Energy Theorem.  In their application,
special solutions of Jang's equation which exhibit blow-up
behavior at apparent horizons played an integral role, and for
some time it has been suggested that these solutions may be
helpful in studying the Penrose Inequality (see [10] for some
problems that can occur with this approach).  For this to be the
case, solutions which blow-up at a given apparent horizon must
always be shown to exist.  In fact such a theorem has recently
been established by Metzger in [11].  More precisely Metzger has
shown that given an initial data set containing an outermost
future (or past) apparent horizon, there exists a smooth solution
of Jang's equation outside of the outermost apparent horizon which
blows-up to $+\infty$ ($-\infty$) in the form of a cylinder over
the horizon, and vanishes at spatial infinity.  Here an outermost
future (past) apparent horizon refers to a future (past) apparent
horizon outside of which there is no other apparent horizon; such
a horizon may have several components. We will denote the Jang
surface associated with the given blow-up solution of Jang's
equation by $(\overline{M},\overline{g})$, and its connection by
$\overline{\nabla}$.  We will show\medskip

\textbf{Theorem 1.2.}  \textit{Let $(M,g,k)$ be an asymptotically
flat initial data set for the Einstein equations satisfying the
dominant energy condition $\mu\geq|J|$.  If the boundary consists
of an outermost future (past) apparent horizon with components of
area $|\partial_{i}M|$, $i=1,\ldots,n$, then the ADM energy
satisfies}
\begin{equation*}
\mathrm{E}(g)\geq\frac{\sigma}{2(1+\sigma)}\sum_{i=1}^{n}\sqrt{\frac{|\partial_{i}
M|}{\pi}}
\end{equation*}
\textit{where}
\begin{equation*}
\sigma=\left(\sum_{i=1}^{n}\sqrt{4\pi|\partial_{i}M|}\right)^{-1}\inf
\parallel\overline{\nabla} v\parallel^{2}_{L^{2}(\overline{M})},
\end{equation*}
\textit{with the infimum taken over all $v\in C^{\infty}(M)$ such
that $v(x)\rightarrow 0$ as $x\rightarrow\partial M$ and
$v(x)\rightarrow 1$ as $|x|\rightarrow\infty$.}\medskip

\textbf{Remark 1.3.}  \textit{Although the hypotheses require a
boundary consisting entirely of future or entirely of past
apparent horizons, our proof gives a bit more.  Namely when both
types are present the same result holds, where
$\{\partial_{i}M\}_{i=1}^{n}$ consists entirely of future or
entirely of past apparent horizons.}\medskip

  An important point to note concerning Theorem 1.2 is that the
case of equality is not considered.  The reasons for this are the
following.  First the Jang equation is designed to embed the
initial data into Minkowski space if equality were to occur (as is
done in the Positive Energy Theorem), and so there is no chance of
obtaining and embedding into the Schwarzchild spacetime in this
situation, as the Penrose Inequality demands.  Moreover, it will
in fact be shown that the case of equality can never be achieved.
This implies that the current result is not optimal (unlike
Theorem 1.1), and suggests that there may be a better choice of
boundary conditions for the Jang equation which does yield an
optimal result.\par
  Another point to note is that the constant $\sigma$ is
dimensionless, and so is actually independent of the area of the
boundary $\partial M$.  Furthermore $\sigma$ never vanishes, and
therefore Theorem 1.2 does give a positive lower bound for the ADM
mass in terms of the area of the apparent horizon, which is
consistent with the spirit of the Penrose Inequality.  Moreover
the theorem may be generalized to the setting of initial data
containing a trapped surface, to give a positive lower bound for
the ADM mass in terms of the least area required to enclose the
trapped surface.  To see this recall that Andersson and Metzger
[1], and Eichmair [4], have shown that the existence of a future
(past) trapped surface in an asymptotically flat initial data set
implies the existence of an outermost future (past) apparent
horizon.  One may then apply Theorem 1.2 to obtain the desired
result.\par
  The proof of Theorem 1.2 closely follows that of Theorem 1.1.
The main difference, or new idea, is to employ blow-up solutions
for Jang's equation in an appropriate way.  However the argument
still relies on the following version of the Positive Energy
Theorem due to Herzlich.\medskip

\textbf{Theorem 1.4 [7].}  \textit{Let $(M,g)$ be a 3-dimensional
asymptotically flat Riemannian manifold with nonnegative scalar
curvature.  If the boundary $\partial M$ consists of $n$
components having spherical topology and mean curvature
(calculated with respect to the normal pointing inside $M$)
satisfying $H_{\partial_{i}M}\leq\sqrt{16\pi/|\partial_{i} M|}$,
$1\leq i\leq n$, then $E(g)\geq 0$ and when equality occurs $g$ is
flat.}\medskip

\textbf{Remark 1.5.}  \textit{The statement of Herzlich's original
theorem only allowed the boundary $\partial M$ to have one
component.  However the same spinor proof may easily be extended
to allow for finitely many components as in Theorem 1.4.}\medskip

  In continuing with the outline of proof for Theorem 1.2, there
are three primary steps.  The first is to deform the initial data
by constructing a blow-up solution of the Jang equation, which as
mentioned above has already been established.  This deformation
yields a positivity property for the scalar curvature of the Jang
metric $\overline{g}$.  The next step entails cutting off the
cylindrical ends of the blown-up Jang surface at a height $T$ to
obtain a manifold with boundary $\overline{M}_{T}$, and then
making a conformal deformation
$(\overline{M}_{T},\widehat{g}_{T}:=u_{T}^{4}\overline{g})$ to
obtain a manifold with zero scalar curvature and with each
boundary component satisfying
$\widehat{H}_{\partial_{i}\overline{M}_{T}}=\sqrt{16\pi/|\partial_{i}
\overline{M}_{T}|_{\widehat{g}_{T}}}$, $1\leq i\leq n$.  Existence
of a conformal factor $u_{T}$ satisfying these properties will be
established by a variational argument, which heavily depends on
the positivity property for the scalar curvature of the Jang
metric as well as the blow-up behavior of the Jang surface at the
horizon.  One may then undertake the last step, which consists of
applying Theorem 1.4 to obtain $E(\widehat{g}_{T})\geq 0$. The
desired lower bound for $E(g)=E(\overline{g})$ is then produced by
estimating the difference $E(\overline{g})-E(\widehat{g}_{T})$ and
letting $T\rightarrow\infty$.
%\par
%  We note in passing that the above difference of masses is proportional
%to the $\frac{1}{r}$-rate at which the conformal factor tends to 1
%as $|x|\rightarrow\infty$ (the conformal factors converge to a
%limit conformal factor as $T\rightarrow\infty$, and the
%$\frac{1}{r}$-rate which we refer to here is that of the limit
%conformal factor).  This number measures the size of the lower
%bound for the ADM mass, or interpreted in another way, it measures
%the contribution to the ADM mass of the energy inside the horizon.
%Thus it may be possible to define a quasi-local mass out of this
%quantity.  As mentioned above the arguments for Theorem 1.2 may be
%generalized to the case in which the boundary consists of trapped
%surfaces. Furthermore, in the time symmetric case these arguments
%may be generalized so that the quasi-local mass applies to a
%boundary satisfying $H\leq\sqrt{16\pi/|\partial M|}$.

\begin{center}
\textbf{2.  The Jang Surface}
\end{center} \setcounter{equation}{0}
\setcounter{section}{2}

   The goal of this section is to give a precise description of
the blow-up solution to the Jang equation, as well as to record
certain qualitative properties of the resulting Jang surface.  Let
us first recall some basic facts.  The Jang surface $\overline{M}$
is given by a graph $t=f(x)$ in the product manifold
$(M\times\mathbb{R},g+dt^{2})$, and so has induced metric
$\overline{g}=g+df^{2}$. The function $f$ is required to satisfy
the Jang equation:
\begin{equation}
\overline{g}^{ij}\left(\frac{\nabla_{ij}f}{\sqrt{1+|\nabla_{g}f|^{2}}}-k_{ij}\right)=0.
\end{equation}
Here $\nabla_{ij}$ denote second covariant derivatives with
respect to $g$ and
\begin{equation*}
\overline{g}^{ij}=g^{ij}-\frac{f^{i}f^{j}}{1+|\nabla_{g}f|^{2}}
\end{equation*}
is the inverse matrix for $\overline{g}_{ij}$ with
$f^{i}=g^{ij}\nabla_{i}f$, and therefore Jang's equation simply
asserts that the mean curvature of the graph is equal to the trace
of $k$ over the graph (assuming that the tensor $k$ has been
extended trivially to all of $M\times\mathbb{R}$).  The motivation
for solving Jang's equation is to obtain a positivity property for
the scalar curvature of the Jang surface.  In particular, if $f$
satisfies equation (2.1) then the scalar curvature of
$\overline{g}$ has the following expression (see [13])
\begin{equation}
\overline{R}=16\pi(\mu-J(w))+
|\overline{h}-k|_{\overline{g}}^{2}+2
|q|_{\overline{g}}^{2}-2\mathrm{div}_{\overline{g}}(q),
\end{equation}
where
\begin{equation*}
w_{i}=\frac{\nabla_{i}f}{\sqrt{1+|\nabla_{g}f|^{2}}},\text{
}\text{ }\text{ }\text{
}q_{i}=\frac{f^{j}}{\sqrt{1+|\nabla_{g}f|^{2}}}(\overline{h}_{ij}-k_{ij}),
\end{equation*}
and $\overline{h}$ is the second fundamental form of
$\overline{M}$. In addition to the positivity property for the
scalar curvature, we will require the Jang surface to exhibit
blow-up behavior at $\partial M$ in order to construct the
conformal factor described in the introduction.  It turns out that
such a solution always exists as long as $\partial M$ is an
outermost horizon.\medskip

\textbf{Theorem 2.1 [11].}  \textit{Suppose that $\partial M$ is
an outermost future (past) apparent horizon.  Then there exists an
open set $\Omega\subset M$ (with $(M-\Omega)\cap\partial
M=\emptyset$) and a smooth function
$f:\Omega\rightarrow\mathbb{R}$ satisfying (2.1), such that
$\partial\Omega-\partial M$ consists of past (future) apparent
horizons, $\overline{M}=\mathrm{graph}(f)$ is asymptotic to the
cylinders $\partial M\times\mathbb{R}_{+}$ ($\partial
M\times\mathbb{R}_{-}$) and $\partial\Omega\times\mathbb{R}_{-}$
($\partial\Omega\times\mathbb{R}_{+}$), and $f(x)\rightarrow 0$ as
$|x|\rightarrow\infty$.}\medskip

  This theorem yields the desired blow-up behavior at $\partial M$
with the added feature that blow-up may occur elsewhere at
$\partial\Omega-\partial M$ as well, if $M$ contains apparent
horizons of the other type (with respect to $\partial M$). However
the hypotheses of Theorem 1.2 do not allow for such extra
horizons, so that in fact $\Omega=M$.  We remark that the sole
reason for prohibiting extra horizons in the initial data is to
ensure that each component of $\partial\Omega$ has spherical
topology, which is needed when applying the Positive Energy
Theorem, Theorem 1.4.  Thus one could allow apparent horizons of
the other type, if they have spherical topology.

  The other goal of this section is to record the decay rate for
certain geometric quantities associated with the Jang surface.
Since the solution of Jang's equation blows-up at $\partial M$ in
the form of a cylinder, in a neighborhood of each boundary
component the Jang surface may be foliated by the level sets
$t=f(x)$, which we denote by $\overline{\Sigma}_{t}$.  Similarly
in a neighborhood of each boundary component, $M$ may be foliated
by the projection of the level sets $\overline{\Sigma}_{t}$ onto
$M$, which we denote by $\Sigma_{t}$. We can then introduce
coordinates $(r,\xi^{2},\xi^{3})$ in such a neighborhood of each
component, where $r=|t|^{-1}$ and $\xi^{2}$, $\xi^{3}$ are
coordinates on a 2-sphere.  Note that as $r\rightarrow 0$ the
projections $\Sigma_{r}$ converge to their associated component of
$\partial M$. Furthermore the $r$-coordinate may be chosen
orthogonal to its level sets, so that the initial data metric
takes the form
\begin{equation*}
g=g_{11}dr^{2}+\sum_{i,j=2}^{3}g_{ij}d\xi^{i}d\xi^{j}.
\end{equation*}

\textbf{Lemma 2.2.}  \textit{Consider the level sets
$\overline{\Sigma}_{r}$ of the blown-up Jang surface near a
component of $\partial M$.  If
$\overline{H}_{\overline{\Sigma}_{r}}$ denotes the mean curvature
of $\overline{\Sigma}_{r}$ with respect to the inward pointing
(towards spatial infinity) normal $\overline{N}$, then
$\overline{H}_{\overline{\Sigma}_{r}}-q(\overline{N})\rightarrow
0$ as $r\rightarrow 0$.}\medskip

\textit{Proof.}  A calculation in [14] (page 10) shows that
\begin{equation*}
\overline{H}_{\overline{\Sigma}_{r}}-q(\overline{N})
=\sqrt{1+|\nabla_{g}f|^{2}}(H_{\Sigma_{r}}\pm\mathrm{Tr}_{\Sigma_{r}}k)
\mp\frac{\mathrm{Tr}_{\Sigma_{r}}k}{|\nabla_{g}f|+\sqrt{1+|\nabla_{g}f|^{2}}},
\end{equation*}
where $H_{\Sigma_{r}}$ is the mean curvature of $\Sigma_{r}$,
$\mathrm{Tr}_{\Sigma_{r}}k$ is the trace over $\Sigma_{r}$, and
$+$ ($-$) is chosen depending on whether the particular component
of $\partial M$ in question is a future (past) horizon
respectively. From this expression we see that it is enough to
show that the first term on the right-hand side approaches zero as
$r\rightarrow 0$. Fortunately this same expression appears in the
Jang equation, and yields the desired result.  To see this, write
the Jang equation in the coordinates $(r,\xi^{2},\xi^{3})$ to
obtain
\begin{equation*}
\frac{g^{11}}{1+g^{11}f_{,r}^{2}}(f_{,rr}-\Gamma_{11}^{1}f_{,r})
-\sum_{i,j=2}^{3}g^{ij}\Gamma_{ij}^{1}f_{,r}=\sqrt{1+g^{11}f_{,r}^{2}}\left(
\frac{g^{11}}{1+g^{11}f_{,r}^{2}}k_{11}+\sum_{i,j=2}^{3}g^{ij}k_{ij}\right),
\end{equation*}
where $f_{,r}$, $f_{,rr}$ are partial derivatives and
$\Gamma_{ij}^{1}$ are Christoffel symbols for $g$ given by
\begin{equation*}
\Gamma_{11}^{1}=\frac{1}{2}g^{11}\partial_{r}g_{11},\text{ }\text{
}\text{ }\text{ }\Gamma_{ij}^{1}=-\sqrt{g^{11}}h_{ij},\text{
}\text{ }\text{ }\text{ }2\leq i,j\leq 3,
\end{equation*}
with $h_{ij}$ denoting the second fundamental form of
$\Sigma_{r}$.  It follows that
\begin{equation*}
\sqrt{1+|\nabla_{g}f|^{2}}(H_{\Sigma_{r}}\pm\mathrm{Tr}_{\Sigma_{r}}k)
=\pm\frac{g^{11}f_{,rr}}{1+g^{11}f_{,r}^{2}}+O\left(\frac{1}{\sqrt{1+|\nabla_{g}f|^{2}}}\right).
\end{equation*}
Lastly we observe that by definition of the coordinate $r$,
$f(r)=\pm r^{-1}$, and therefore
\begin{equation*}
\sqrt{1+|\nabla_{g}f|^{2}}(H_{\Sigma_{r}}\pm\mathrm{Tr}_{\Sigma_{r}}k)
=O(r)\text{ }\text{ }\text{ as }\text{ }\text{ }r\rightarrow 0.
\end{equation*}
Q.E.D.\medskip

\textbf{Lemma 2.3.}  \textit{The solution of Jang's equation
satisfies the following fall-off condition at spatial infinity:}
\begin{equation*}
|\nabla^{l}f|(x)=O(|x|^{-\frac{1}{2}-l})\text{ }\text{ }\textit{
as }\text{ }\text{ }|x|\rightarrow\infty,\text{ }\text{ }\text{
}\text{ }l=0,1,2.
\end{equation*}
\textit{In particular, the energy of the Jang metric
$\overline{g}$ equals the energy of $g$.}\medskip

\textit{Proof.}  See Schoen and Yau [13].  Q.E.D.

\begin{center}
\textbf{3.  The Conformal Factor}
\end{center} \setcounter{equation}{0}
\setcounter{section}{3}

  In this section we will complete the last preliminary step
before application of the Positive Energy Theorem.  Namely we will
conformally deform the Jang metric to zero scalar curvature on a
portion of the Jang surface, while at the same time prescribing
the mean curvature of its boundary.  The region to be considered
consists of the portion of the Jang surface lying between the
horizontal planes $t=\pm T$, and will be denoted by
$\overline{M}_{T}$.  We then search for a conformal factor $u_{T}$
satisfying the following boundary value problem
\begin{equation}
\overline{\Delta}u_{T}-\frac{1}{8}\overline{R}u_{T}=0\text{
}\text{ }\text{ on }\text{ }\text{ }\overline{M}_{T},
\end{equation}
\begin{equation*}
\partial_{\overline{\nu}}u_{T}+\frac{1}{4}\overline{H}_{\partial_{i}\overline{M}_{T}}u_{T}=\frac{1}{4}
\sqrt{\frac{16\pi}{|\partial_{i}\overline{M}_{T}|_{\widehat{g}_{T}}}}u_{T}^{3}\text{
}\text{ }\text{ on each $i$th component of }\text{ }\text{
}\partial\overline{M}_{T},
\end{equation*}
\begin{equation*}
u_{T}=1+\frac{A_{T}}{|x|}+O(|x|^{-2})\text{ }\text{ }\text{ as
}\text{ }\text{ }|x|\rightarrow\infty,
\end{equation*}
where $A_{T}$ is a constant and
$\overline{H}_{\partial_{i}\overline{M}_{T}}$ denotes mean
curvature with respect to the unit inward normal $\overline{\nu}$
(pointing inside $\overline{M}_{T}$). This ensures that
$(\overline{M}_{T},\widehat{g}_{T}:=u_{T}^{4}\overline{g})$ has
zero scalar curvature $\widehat{R}\equiv 0$ and mean curvature on
each $i$th component of $\partial\overline{M}_{T}$ given by
$\widehat{H}_{\partial_{i}\overline{M}_{T}}
=\sqrt{16\pi/|\partial_{i}\overline{M}_{T}|_{\widehat{g}_{T}}}$.
These two properties, combined with the fact that each component
of $\partial M$ must have spherical topology ([5], [6]), then
guarantee that Theorem 1.4 is applicable.\par
  We shall use a variational argument, just as in [7], to
construct $u_{T}:=1+v_{T}$.  In this regard observe that boundary
value problem (3.1) arises as the Euler-Lagrange equation for the
functional
\begin{equation*}
Q(v)=\frac{1}{2}\int_{\overline{M}_{T}}\left(|\overline{\nabla}v|^{2}+\frac{1}{8}\overline{R}
(1+v)^{2}\right)+\frac{\sqrt{\pi}}{2}\left(\int_{\partial\overline{M}_{T}}(1+v)^{4}\right)^{1/2}
-\frac{1}{8}\int_{\partial\overline{M}_{T}}\overline{H}_{\partial\overline{M}_{T}}(1+v)^{2}.
\end{equation*}
More precisely, we will search for a global minimum over the
weighted Sobolev space
\begin{equation*}
W^{1,2}_{-1}(\overline{M}_{T})=\{v\in
W^{1,2}_{loc}(\overline{M}_{T})\mid
|x|^{l-1}\overline{\nabla}^{l}v\in L^{2}(\overline{M}_{T}),\text{
}\text{ }l=0,1\}.
\end{equation*}

\textbf{Theorem 3.1.}  \textit{Given $T>0$ sufficiently large,
there exists a function $v_{T}\in
W^{1,2}_{-1}(\overline{M}_{T})\cap C^{\infty}(\overline{M}_{T})$
at which $Q$ attains a global minimum.  Moreover $u_{T}=1+v_{T}$
never vanishes and satisfies the asymptotic behavior in
(3.1).}\medskip

\textit{Proof.}  In order to establish the existence (as well as
the regularity and asymptotic behavior) portion of this theorem it
is enough, by the arguments of [7], to show that for $T$
sufficiently large the functional $Q$ is nonnegative. To see this
use formula (2.2) and integrate the divergence term by parts to
find that for any $v\in W^{1,2}_{-1}(\overline{M}_{T})$,
\begin{eqnarray}
Q(v)&\geq&\int_{\overline{M}_{T}}\left(\frac{3}{8}|\overline{\nabla}v|^{2}+\pi(\mu-|J|)(1+v)^{2}\right)+
\frac{\sqrt{\pi}}{2}\left(\int_{\partial\overline{M}_{T}}(1+v)^{4}\right)^{1/2}\\
& &
-\frac{1}{8}\int_{\partial\overline{M}_{T}}(\overline{H}_{\partial\overline{M}_{T}}-q(\overline{N}))(1+v)^{2}.\nonumber
\end{eqnarray}
By Lemma 2.2
$\overline{H}_{\partial\overline{M}_{T}}-q(\overline{N})=O(T^{-1})$,
and a calculation shows that the area of
$\partial\overline{M}_{T}$ agrees with the area of
$\Sigma_{T}\subset M$ which remains bounded as
$T\rightarrow\infty$.  It then follows from Jensen's Inequality
\begin{equation*}
\left(\int_{\partial\overline{M}_{T}}(1+v)^{2}\right)^{2}\leq
|\partial\overline{M}_{T}|\int_{\partial\overline{M}_{T}}(1+v)^{4},
\end{equation*}
that for $T$ sufficiently large $Q$ is nonnegative.\par
  It remains to show that $u_{T}=1+v_{T}$ is strictly positive.
So suppose that $u_{T}$ is not positive and let $D_{-}$ be the
domain on which $u_{T}<0$.  Since $u_{T}\rightarrow 1$ as
$|x|\rightarrow\infty$, the closure of $D_{-}$ must be compact.
Now multiply equation (3.1) through by $u_{T}$ and integrate by
parts to obtain
\begin{equation*}
\int_{D_{-}}|\overline{\nabla}u_{T}|^{2}\leq 0.
\end{equation*}
Note that if $D_{-}\cap\partial\overline{M}_{T}\neq\emptyset$,
then the same arguments used above to show that $Q$ is
nonnegative, must be employed.  It follows that $u_{T}\geq 0$.  To
show that $u_{T}>0$, one need only apply Hopf's Maximum Principle
(the boundary condition of (3.1) must be used to obtain this
conclusion at $\partial\overline{M}_{T}$). Q.E.D.

\begin{center}
\textbf{4.  Proof of Theorem 1.2}
\end{center} \setcounter{equation}{0}
\setcounter{section}{4}

   Here we shall carry out the last step in the proof of Theorem
1.2, namely to apply the Positive Energy Theorem and to compare
the two energies $\mathrm{E}(g)$ and
$\mathrm{E}(\widehat{g}_{T})$. Observe that all the hypotheses of
Theorem 1.4 are satisfied by $(\overline{M}_{T},\widehat{g}_{T})$
so that $\mathrm{E}(\widehat{g}_{T})\geq 0$. Therefore a
straightforward calculation yields
\begin{equation}
\mathrm{E}(g)\geq \mathrm{E}(g)-\mathrm{E}(\widehat{g}_{T})
=\frac{1}{2\pi}\lim_{r\rightarrow\infty}\int_{|x|=r}\partial_{\overline{\nu}}u_{T}.
\end{equation}
Furthermore upon integrating by parts and using boundary value
problem (3.1) we obtain
\begin{equation}
\lim_{r\rightarrow\infty}\int_{|x|=r}\partial_{\overline{\nu}}u_{T}=
\lim_{r\rightarrow\infty}\int_{|x|=r}u_{T}\partial_{\overline{\nu}}u_{T}=2Q(v_{T}).
\end{equation}\par
  Now suppose that
$Q(v_{T})\leq\eta\sum_{i=1}^{n}\sqrt{\pi|\partial_{i}\overline{M}_{T}|}$
for some positive constant $\eta$, where $n$ denotes the number of
components comprising $\partial M$.  Then integrating by parts,
and using arguments such as those found in the proof of Theorem
3.1, shows that there exists a constant $C>0$ independent of $T$
such that
\begin{equation*}
\frac{3}{8}\int_{\overline{M}_{T}}|\overline{\nabla}v_{T}|^{2}+\left(\frac{1-CT^{-1}}{2}\right)
\sum_{i=1}^{n}\sqrt{\frac{\pi} {|\partial_{i}\overline{M}_{T}|}}
\int_{\partial_{i}\overline{M}_{T}}(1+v_{T})^{2}
\leq\eta\sum_{i=1}^{n}\sqrt{\pi|\partial_{i}\overline{M}_{T}|}.
\end{equation*}
However by Young's Inequality
\begin{equation*}
(1+v_{T})^{2}\geq 1-\frac{1}{\delta}+(1-\delta)v_{T}^{2}
\end{equation*}
for any $\delta>0$, and therefore
\begin{eqnarray*}
&
&\frac{3}{8}\int_{\overline{M}_{T}}|\overline{\nabla}v_{T}|^{2}+(1-\delta)
\left(\frac{1-CT^{-1}}{2}\right)\sum_{i=1}^{n}\sqrt{\frac{\pi}
{|\partial_{i}\overline{M}_{T}|}}\int_{\partial_{i}\overline{M}_{T}}v_{T}^{2}\\
&\leq& (\eta-\frac{1}{2}(1-\delta^{-1})(1-CT^{-1}))
\sum_{i=1}^{n}\sqrt{\pi|\partial_{i}\overline{M}_{T}|}.
\end{eqnarray*}
It follows that the left-hand side is nonnegative if
$\delta-1\leq\sigma_{T}$ where
\begin{equation*}
\sigma_{T}=
\frac{\int_{\overline{M}_{T}}|\overline{\nabla}v_{T}|^{2}}{2(1-CT^{-1})
\sum_{i=1}^{n}\sqrt{\frac{\pi}{|\partial_{i}\overline{M}_{T}|}}
\int_{\partial_{i}\overline{M}_{T}}v_{T}^{2}},
\end{equation*}
so that $\eta\geq\delta^{-1}(\delta-1)(1-CT^{-1})/2$ for all such
$\delta$.  In particular by choosing $\delta=1+\sigma_{T}$ we
conclude that
\begin{equation}
Q(v_{T})\geq\frac{\sigma_{T}(1-CT^{-1})}{2(1+\sigma_{T})}\sum_{i=1}^{n}\sqrt{\pi
|\partial_{i}\overline{M}_{T}|}.
\end{equation}
Furthermore combining (4.1), (4.2), and (4.3) produces
\begin{equation}
\mathrm{E}(g)\geq \frac{\sigma_{T}(1-CT^{-1})}{2(1+\sigma_{T})}
\sum_{i=1}^{n}\sqrt{\frac{|\partial_{i}\overline{M}_{T}|}{\pi}}.
\end{equation}\par
  The desired inequality of Theorem 1.2 may be obtained from (4.4)
by letting $T\rightarrow\infty$.  To see this we observe that
(4.1), (4.2), and (3.2) together show that the sequence of
functions $\{u_{T}\}$ is uniformly bounded in
$W^{1,2}_{loc}(\overline{M})$. Thus with the help of elliptic
estimates and Sobolev embeddings, this sequence converges on
compact subsets to a smooth uniformly bounded solution
$u_{\infty}$ of
\begin{equation*}
\overline{\Delta}u_{\infty}-\frac{1}{8}\overline{R}u_{\infty}=0\text{
}\text{ }\text{ on }\text{ }\text{ }\overline{M},\text{ }\text{
}\text{ }\text{
}u_{\infty}=1+\frac{A_{\infty}}{|x|}+O(|x|^{-2})\text{ }\text{
}\text{ as }\text{ }\text{ }|x|\rightarrow\infty.
\end{equation*}
However since $\overline{M}$ approximates a cylinder on regions
where it blows-up, comparison with a bounded solution of the same
equation on the cylinder (as is done in [13]) shows that
$u_{\infty}(x)\rightarrow 0$ as $x\rightarrow\partial M$; in fact
the decay rate is of exponential strength.  Therefore (with a bit
more effort) $\sigma_{T}\rightarrow\sigma_{\infty}\geq \sigma$ and
$|\partial_{i}\overline{M}_{T}|\rightarrow|\partial_{i}M|$, $1\leq
i\leq n$, as $T\rightarrow\infty$.  This completes the proof of
Theorem 1.2.\par
  Lastly we analyze what happens when equality occurs in Theorem
1.2.  By slightly modifying the arguments of this section in this
special case, we find that
\begin{equation*}
\int_{\overline{M}}|\overline{\nabla}u_{\infty}|^{2}=0,
\end{equation*}
and therefore $u_{\infty}$ must be constant.  However this is
impossible since
\begin{equation*}
u_{\infty}(x)\rightarrow\begin{cases}
1 & \text{as $|x|\rightarrow\infty$},\\
0 & \text{as $x\rightarrow\partial M$}.
\end{cases}
\end{equation*}
We conclude that the case of equality cannot occur.

\begin{center}
\textbf{References}
\end{center}

\noindent 1.\hspace{.07in} L. Andersson, and J. Metzger,
\textit{The area of horizons and the trapped region},
\par\hspace{-.01in} preprint, arXiv:0708.4252, 2007.\medskip

\noindent 2.\hspace{.07in} H. Bray, \textit{Proof of the
Riemannian Penrose conjecture using the positive mass}
\par\hspace{-.01in} \textit{theorem,} J. Differential Geom.,
$\mathbf{59}$ (2001), 177-267.\medskip

\noindent 3.\hspace{.07in} H. Bray, and M. Khuri, \textit{PDE's
which imply the Penrose conjecture}, in preparation,
\par\hspace{-.01in} 2008.\medskip

\noindent 4.\hspace{.07in} M. Eichmair, \textit{The plateau
problem for apparent horizons},  preprint, arXiv:0711.4139,
\par\hspace{-.01in} 2007.\medskip

\noindent 5.\hspace{.07in} G. Galloway, \textit{Rigidity of
marginally trapped surfaces and the topology of black}
\par\hspace{-.01in} \textit{holes}, Commun. Anal. Geom., $\mathbf{16}$
(2008), 217-229.\medskip

\noindent 6.\hspace{.07in} G. Galloway, and R. Schoen, \textit{A
generalization of Hawking's black hole topology}
\par\hspace{-.01in} \textit{theorem to higher dimensions},
Commun. Math. Phys., $\mathbf{266}$ (2006), no. 2, 571-
\par\hspace{-.01in} 576.\medskip

\noindent 7.\hspace{.07in} M. Herzlich, \textit{A Penrose-like
inequality for the mass of Riemannian asymptotically}
\par\hspace{-.01in} \textit{flat manifolds}, Commun. Math. Phys.,
$\mathbf{188}$ (1997), 121-133.\medskip

\noindent 8.\hspace{.07in} G. Huisken, and T. Ilmanen, \textit{The
inverse mean curvature flow and the Rieman-}
\par\hspace{-.01in} \textit{nian Penrose inequality}, J. Differential Geom.,
$\mathbf{59}$ (2001), 353-437.\medskip

\noindent 9.\hspace{.07in} P.-S. Jang, \textit{On the positivity
of energy in General Relativity}, J. Math. Phys., $\mathbf{19}$
\par\hspace{-.01in} (1978), 1152-1155.\medskip

\noindent 10.  E. Malec, and N. \'{O} Murchadha, \textit{The Jang
equation, apparent horizons, and the}
\par\hspace{-.01in} \textit{Penrose inequality}, Class. Q. Grav.,
$\mathbf{21}$ (2004), 5777-5787.\medskip

\noindent 11.  J. Metzger, \textit{Blowup of Jang's equation at
outermost marginally trapped surfaces},
\par\hspace{-.01in} preprint, arXiv:0711.4753, 2008.\medskip

\noindent 12.  R. Penrose, \textit{Naked singularities}, Ann. N.
Y. Acad. Sci., $\mathbf{224}$ (1973), 125-134.\medskip

\noindent 13.  R. Schoen, and S.-T. Yau, \textit{Proof of the
positive mass theorem II}, Commun. Math.
\par\hspace{-.01in} Phys., $\mathbf{79}$ (1981), no. 2, 231-260.\medskip

\noindent 14. S.-T. Yau, \textit{Geometry of three manifolds and
existence of black hole due to} \par\hspace{-.01in}
\textit{boundary effect}, Adv. Theor. Math. Phys., $\mathbf{5}$
(2001), no. 4, 755-767.

\begin{center}
\textbf{Erratum}
\end{center} %\setcounter{equation}{0}

Two errors in %\cite{Khuri}
this paper have been pointed out by Jan Metzger. The first concerns the constant
$\sigma$ in the statement of Theorem 1.2. Namely, this constant is in fact zero. To see this fix large $L>0$. Let $z_{L}\in C^{\infty}
(\overline{M})$ be identically 1 up to height $L$ (in
$M\times\mathbb{R}$), and be zero above height $2L$, with
$|\overline{\nabla}z_{L}|< 2L^{-1}$.  It follows that
$\parallel\overline{\nabla}z_{L}\parallel_{L^{2}(\overline{M})} \sim
L^{-1/2}$. Hence, we find that $\sigma=0$ by letting
$L\rightarrow\infty$. The second error concerns Lemma 2.2. Namely, the quantity $\overline{H}-q(\overline{N})$ may not necessarily approach
zero as $r\rightarrow 0$. This may be seen by examining certain spherically symmetric examples.  Lemma 2.2 is not necessary in order to establish
the main result. More precisely, the Penrose-like inequality with a corrected definition of $\sigma$ may be derived
from the ideas in this paper, together with some additional arguments. The details will be given by the author in a separate paper (arXiv:1308.3591).

%\begin{thebibliography}{99}

%\bibitem{Khuri} M. A. Khuri, \textit{A Penrose-like
%inequality for general initial data sets}, Comm. Math. Phys.,
%\textbf{290}(2009), no. 2, 779-788.

%\end{thebibliography}

\bigskip\bigskip\footnotesize

\noindent\textsc{Department of Mathematics, Stony Brook
University, Stony Brook, NY 11794}\par

\noindent\textit{E-mail address}: \verb"khuri@math.sunysb.edu"

\end{document}